\documentclass[11pt,twoside]{amsart}
\usepackage{amsmath, amsthm, amscd, amsfonts, amssymb, graphicx, color}
\usepackage[bookmarksnumbered, plainpages]{hyperref}
\usepackage{color}
\numberwithin{equation}{section}

\begin{document}

\begin{center}

{\Large{\textbf{\textsc{On the Preconditioned AOR Iterative Method for Z-Matrices}}}} \vspace{0.2cm}
\\
\vspace{0.5cm}

{ Davod Khojasteh Salkuyeh$^\dag$\footnote{Corresponding author}, Mohsen Hasani$^\ddag$ and
Fatemeh Panjeh Ali Beik$^\S$} \\[0.2cm]

$^\dag${\small \textit{Faculty of Mathematical Sciences, University of Guilan}\\\textit{P.O. Box 1914, Rasht, Iran}\\
\texttt{email:~salkuyeh@gmail.com, khojasteh@guilan.ac.ir}} \\[0.2cm]

$^\ddag${\small \textit{Faculty of Science, Department of Mathematics, Islamic Azad University,}\\\textit{Shahrood, Iran}\\
\texttt{email:~hasani.mo@gmail.com}} \\ [0.2cm]

$^\S${\small \textit{Department of Mathematics, Vali-e-Asr University of Rafsanjan,}\\\textit{Rafsanjan, Iran}\\
\texttt{email:~f.beik@vru.ac.ir}} \\

\end{center}

\bigskip

\medskip

\begin{center}
 {\bf Abstract}
\end{center}

\indent Consider a general class of preconditioners which are nonsingular, nonnegative
and has unit diagonal entries. In the various recently published papers, the authors have applied particular
preconditioners from this class of preconditioners  to propose the preconditioned AOR
methods for solving a linear system of equations with a unit Z-matrix coefficient matrix.
The main purpose of this paper is to present a comparison result among the preconditioned
AOR methods where the preconditioners are more general and pertain to the delineated class
of preconditioners. Numerical experiments for corresponding preconditioned GMRES methods are
reported to illustrate the theoretical results.

\bigskip
\noindent\textit{AMS Subject Classification :} 65F10, 65F50.\\
\noindent\textit{Keywords}:  Linear system of equations, Preconditioner, AOR iterative method, Z-matrix, Comparison result.

\pagestyle{myheadings}
\markboth{\rightline {\scriptsize  D. K. Salkyeh, M. Hasani, and F. P. A. Beik}}
         {\leftline{\scriptsize On the preconditioned AOR iterative method for Z-matrices
}}

\bigskip

\bigskip

\bigskip

\section{Introduction}
\label{SEC1}

We consider the following linear system of  equations
\begin{equation}\label{eq001}
    Ax=b,\\
\end{equation}
where $A=(a_{ij})\in \mathbb{R}^{n\times n}$ is nonsingular and $b\in \mathbb{R}^n$. By means of the splitting $A=M-N$ in which  $M,N\in \mathbb{R}^{n\times n}$ and $M$ is nonsingular, a general stationary iterative method for solving Eq. (\ref{eq001}) is expressed as follows:
\begin{equation}\label{eq0001}
x^{(k+1)}=M^{-1}Nx^{(k)}+M^{-1}b,\quad k=0,1,2,\ldots,
\end{equation}
where the initial vector $x^{(0)}$  is given and ${\mathcal{L}}=M^{-1}N$ is called the iteration matrix.

It is well-known that the iterative method \eqref{eq0001} is convergent for each arbitrary choice of the starting vector $x^{(0)}$  if and only if $\rho(\mathcal{L})<1$. Here, the notation $\rho(X)$ represents the spectral radius of the matrix $X$. In this paper, we assume that $a_{ii}\neq 0$ for $i=1,2,\ldots,n$. Therefore, without loss of generality, we may presume that all of the diagonal entries of the coefficient matrix $A$ are equal to one. In this situation, we split the matrix $A$ into the subsequent  form
\begin{equation}\label{eq002}
A=I-L-U,
\end{equation}
where $I$ stands for the  identity matrix, $-L$ and $-U$ are strictly lower and strictly upper triangular matrices, respectively. The accelerated overrelaxation (AOR) iterative method for solving  Eq. (\ref{eq001}) is specified  by (for further details see \cite{Hadji,Song1})
\[
x^{(k+1)}=\mathcal{L}_{\gamma,\omega}x^{(k)}+\omega(I-\gamma L)^{-1}b,
\]
in which
\[
\mathcal{L}_{\gamma,\omega}=(I-\gamma L)^{-1}[(1-\omega)I+(\omega-\gamma)L+\omega U],
\]
where $\omega$ and $\gamma$ are real parameters and $\omega\neq 0$. The AOR iterative method incorporates  the Jacobi, Gauss-Seidel and the SOR iterative methods as special cases for  certain values of the parameters $\omega$ and $\gamma$, see \cite{Hadji}.

In order to ameliorate the convergence rate of an iterative method, one may apply it to the  preconditioned linear system $PAx=Pb$. Here, the matrix $P$ is called a preconditioner. In the literature, the application of the several kinds of preconditioners have been investigated widely for the stationary iterative methods.
Nevertheless, many of these preconditioners are the special cases of a  general class of the  preconditioners which has been recently examined by Wang and Song \cite{Wang}, (see for example \cite{Evans,Huang1,Kotakemori,Li1,Li2,Li3,Wang,Zheng}). In \cite{Wang}, the authors have proposed a general preconditioner $P$ which is nonsingular, nonnegative and has unit diagonal entries. More precisely, the authors have investigated the properties of the preconditioners of the form
\begin{equation}\label{eq003}
P=(p_{ij})=(-\alpha_{ij} a_{ij}),
\end{equation}
where  $p_{ii}=1$ and $0\leq \alpha_{ij} \leq 1$ for $i,j=1,\ldots,n$ ($i\neq j$).

In \cite{Noutsos}, Noutsos and Tzoumas have studied the performance  of a  family of the preconditioners which have the following form
\[
P_1  = \left( {\begin{array}{*{20}c}
   1 & {} & {} & { - a_{1k_1} } & {} & {} & {}  \\
   {} & 1 & {} & {} & {} & { - a_{2k_2 } } & {}  \\
   {} & {} & 1 & {} & { - a_{3k_3 } } & {} & {}  \\
   {} & {} & {} &  \ddots  & {} & {} & {}  \\
   {} & {} & {} & {} &  \ddots  & {} & {}  \\
   {} & {} & {} & {} & {} &  \ddots  & {}  \\
   {} & {} & {} & { - a_{nk_n } } & {} & {} & 1  \\
\end{array}} \right),
\]
where $k_i\in \left\{ {1,2,\ldots,i - 1,i + 1,\ldots,n} \right\},$ $i = 1,2,\ldots,n.$
\noindent In order to improve the rate of convergence of the AOR iterative method, Wang et al. \cite{WangG} have extended the above preconditioners to a new one. More precisely,  the authors have handled the following preconditioner
\[
P' = \left( {\begin{array}{*{20}c}
   1 & { - a_{12} } &  \cdots  & { - a_{1,n - 1} } & { - a_{1n} }  \\
   { - a_{21} } & 1 &  \cdots  & { - a_{2,n - 1} } & { - a_{2n} }  \\
    \vdots  &  \vdots  &  \vdots  &  \vdots  &  \vdots   \\
   { - a_{n1} } & { - a_{n2} } &  \cdots  & { - a_{n,n - 1} } & 1  \\
\end{array}} \right).
\]
It has been elaborated that the preconditioned AOR method with the precondiotioner $P'$ outperforms the preconditioned AOR method with the precondiotioner $P_1$. Here, we would like to point out  that our mentioned preconditioner incorporates $P'$. As a matter of fact, the preconditioner $P$ reduces to $P'$  when $\alpha_{ij}=1$ for $i,j=1,2,\ldots,n$.

In this paper, it is  both theoretically and experimentally shown that under some certain assumptions, among the preconditioners of the form (\ref{eq003}), the preconditioner obtained by setting  $\alpha_{ij}=1$ for $i,j=1,2,\ldots,n$ ($i\neq j$) surpasses the other preconditioners.

Before ending this section, we first present some notations. Afterwards, some useful definitions and preliminaries are recollected.

A matrix $A=(a_{ij})\in \mathbb{R}^{n\times n}$ is said to be nonnegative and denoted by $A\geq 0$ if
$a_{ij}\geq 0$ for  $i,j=1,2,\ldots,n$. A matrix $A$ is called positive
and represented by $A\gg 0$ if all of its entries are positive. If $A\geq 0$, then  the well-known Perron-Frobenius theorem implies that $\rho(A)$ is an eigenvalue of $A$, see \cite{Berman}. In addition, corresponding to $\rho(A)$, the matrix $A$ has a nonnegative eigenvector called a Perron vector of $A$.

In the following, we state some definitions and theorems which are utilized throughout of the paper.

\vspace{0.2cm}

\noindent\textbf{Definition 1.1.} A matrix $A=(a_{ij})\in \mathbb{R}^{n\times n}$ is an \textit{Z-matrix} if $a_{ij}\leq 0$ for $i\neq j$.

\bigskip

\noindent\textbf{Definition 1.2.} A Z-matrix $A$ is said to be  an \textit{M-matrix} if $A$ is nonsingular and $A^{-1}\geq 0$.

\bigskip

\noindent\textbf{Definition 1.3.} Suppose that the matrix $A\in \mathbb{R}^{n \times n}$ is given. The representation $A=M-N$ is called a \textit{splitting} of $A$ if $M$ is nonsingular. The splitting $A=M-N$ is called\\
(a) \textit{convergent} if $\rho(M^{-1}N)<1$;\\
(b) \textit{weak regular} if $M^{-1}\geq 0$ and $M^{-1}N\geq 0$;\\
(c) an \textit{M-splitting} of $A$ if $M$ is an M-matrix and $N\geq 0$.

\bigskip

\noindent\textbf{Definition 1.4.} A real matrix $A$ is called \textit{monotone }if $Ax\ge 0$ implies $x\geq 0$.

\bigskip

\noindent\textbf{Lemma 1.1.} \cite[Lemma 3.2 ]{Li3} \textit{Let $A=M-N$ be an M-splitting of $A$. Then $\rho(M^{-1}N)<1$ if and only if $A$ is an M-matrix.}

\bigskip

\noindent\textbf{Lemma 1.2.} \cite[Lemma 1.6 ]{Wu} \textit{Let $A$ be a Z-matrix. Then, $A$ is an M-matrix if and only if there is a positive vector $x$ such that $Ax\gg 0$}.

\vspace{0.2cm}

The following lemma can be  instantly deduced  from the theoretical results proved in \cite{Woznicki}.

\vspace{0.2cm}

\noindent\textbf{Lemma 1.3.}  \textit{ Let $A=M_1-N_1=M_2-N_2$ be two convergent weak regular splittings
of $A$ where $A^{-1}\ge (>) 0$ , if $M_1^{-1} \ge (>) M_2^{-1}$  then
$\rho (M_1^{-1}N_1)\le (<) \rho (M_2^{-1}N_2).$}

\bigskip

\noindent\textbf{Definition 1.5.}   A matrix $A$ is said to be reducible if there is a
permutation matrix $P$ such that $PAP^T$ is a block upper
triangular matrix. Otherwise, it is irreducible.

\vspace{0.2cm}

The following lemma has been originally established by Varga \cite{Varga}. It  provides an easier way to check whether a matrix is irreducible or not.

\vspace{0.2cm}

\noindent\textbf{Lemma 1.4.} \label{lir} A matrix $A$ is irreducible if the directed graph associated to $A$ is strongly connected.

\vspace{0.2cm}

\noindent\textbf{Note.} Throughout this paper the directed graph of matrix $A$ is denoted by $\mathcal{G}(A)$.

\vspace{0.2cm}

The reminder of this paper is organized as follows. In Section 2, we concentrate on a general class of preconditioners to speed up the convergence rate of the AOR method for solving Z-matrix linear system of equations. It is theoretically demonstrated that under certain conditions in the mentioned class of preconditioners one of the preconditioner outperforms other preconditioners. In Section 3, some numerical results are presented  which confirm our theoretical results. Finally, the paper is finished with a succinct conclusion in Section 4.

\section{Main results}
\label{SEC2}

Let us consider the linear system of equations (\ref{eq001}) in which the coefficient matrix $A=(a_{ij})\in \mathbb{R}^{n\times n}$ is given such that $a_{ii}=1$ for $i=1,2,\ldots,n$. In this section, we examine the influence of the preconditioner $\tilde{P}=(\tilde{p}_{ij})\in \mathbb{R}^{n\times n}$ on Eq. (\ref{eq001}) with
\[
\tilde{p}_{ij}=\left\{
  \begin{array}{ll}
    -\alpha_{ij}a_{ij}, & \textrm{if} ~ i \neq j, \\
    \quad 1  , & \textrm{otherwise},
  \end{array}
\right.
\]
where $\alpha_{ij}\in \mathbb{R}$ for $i\neq j$. Let us split the preconditioner $\tilde{P}$ into $\tilde{P}=I+L(\alpha)+U(\alpha)$, in which $I$ is the identity matrix and $L(\alpha)$ and $U(\alpha)$ are strictly lower and strictly upper triangular matrices, respectively. Presume that $\tilde{A}=\tilde{P}A=(I+L(\alpha)+U(\alpha))A$ and
\begin{eqnarray*}
  L(\alpha)U \hspace{-0.2cm} &=& \hspace{-0.2cm} G_1(\alpha)+ E_1(\alpha) +F_1(\alpha), \\
  U(\alpha)L \hspace{-0.2cm} &=& \hspace{-0.2cm} G_2(\alpha)+ E_2(\alpha) +F_2(\alpha),
\end{eqnarray*}
where $E_1(\alpha)$ and
$E_2(\alpha)$ are diagonal matrices, $F_1(\alpha)$ and
$F_2(\alpha)$ are strictly lower triangular matrices and
$G_1(\alpha)$ and $G_2(\alpha)$
are strictly upper triangular matrices.

In this case, the matrix $\tilde{A}$ can be decomposed as $\tilde{A}=\tilde{D} -\tilde{L}-\tilde{U}$. Here, the matrices $\tilde{D}$, $\tilde{L}$ and $\tilde{U}$ are respectively diagonal, strictly lower and strictly upper triangular matrices and defined by
\begin{eqnarray*}
  \tilde{D} \hspace{-0.2cm} &=& \hspace{-0.2cm} I- E_1(\alpha)- E_2(\alpha), \\
  \tilde{L} \hspace{-0.2cm} &=& \hspace{-0.2cm} L-L(\alpha)+ L(\alpha)L+ F_1(\alpha)+ F_2(\alpha),  \\
  \tilde{U} \hspace{-0.2cm} &=& \hspace{-0.2cm} U+ G_1(\alpha)- U(\alpha)+ G_2(\alpha)+ U(\alpha)U.
\end{eqnarray*}
If the matrix $\tilde{D}-\gamma\tilde{L}$ is nonsingular, the AOR iteration matrix for solving  the preconditioned system $\tilde{P}Ax=\tilde{P}b$ can be written as
\[
\tilde{\mathcal{L}}_{\gamma,\omega} =(\tilde{D}-\gamma\tilde{L} )^{-1}[(1-\omega)\tilde{D}+
(\omega-\gamma)\tilde{L}+\omega\tilde{U}].
\]

\bigskip

\noindent\textbf{Theorem  2.1.} \textit{Let $A$ be a Z-matrix and $\alpha_{ij}\in[0,1]$ for $1\leq i\neq j\leq n$. Then, $A$ is an M-matrix if and only if $\tilde{A}$ is and M-matrix.}

\bigskip

\noindent\textbf{Proof.} We may prove the theorem in an analogous manner employed in the proof of Lemma 3.3 in \cite{Li3}. Let $A$ be an M-matrix and $\tilde{A}=\tilde{P}A=(\tilde{a}_{ij})$. Straightforward computations reveal that
\begin{equation}\label{eq004}
\tilde{a}_{ij}=\left\{
  \begin{array}{ll}
    1-\displaystyle\sum_{k=1,k\neq i}^{n} \alpha_{ik}a_{ik}a_{ki}, & 1\leq i=j\leq n, \\
    a_{ij}-\displaystyle\sum_{k=1,k\neq i}^{n} \alpha_{ik}a_{ik}a_{kj}, & 1\leq i\neq j\leq n. \\
  \end{array}
\right.
\end{equation}
By the assumption $A$ is an  Z-matrix, thus  $\tilde{a}_{ij}\leq 0$ for $i, j=1,2,\ldots,n$ and $i\neq j$. Therefore, we conclude that $\tilde{A}$ is also an Z-matrix. By Lemma 1.2, there exists a positive vector $x$ such that $Ax\gg 0$. On the other hand, we have $\tilde{A}=(I+L(\alpha)+U(\alpha))Ax\gg 0$. Invoking Lemma 1.2, we deduce that that $\tilde{A}$ is also an M-matrix.

Conversely, let $\tilde{A}$ be an M-matrix. Evidently, $\tilde{A}^T$ is also an M-matrix and Lemma 1.2 implies the existence of a positive vector $x$ for which  $\tilde{A}^Tx \gg 0$, i.e., $A^T(I+L(\alpha)^T+U(\alpha)^T)x \gg 0$. For simplicity, we set $y=(I+L(\alpha)^T+U(\alpha)^T)x$. It is not difficult to see that $y\gg 0$. Thus,  Lemma 1.2 indicates that $A^T$ is an M-matrix. As a result, $A$ is an M-matrix which completes the proof. \quad $\Box$

\bigskip

\noindent\textbf{Theorem  2.2.} \textit{Let $A=(a_{ij})\in \mathbb{R}^{n\times n}$ be a nonsingular Z-matrix, $0\leq\gamma\leq\omega\leq 1$,  $\omega\neq 0$ and $\alpha_{ij} \in [0,1]$ for $1\leq i \neq j \leq n$.
If $\rho(\mathcal{L}_{\gamma,\omega})<1$, then $\rho(\tilde{\mathcal{L}}_{\gamma,\omega})\leq \rho(\mathcal{L}_{\gamma,\omega})<1$.}

\bigskip

\noindent\textbf{Proof.} Under the assumptions of the theorem, it is easy to see that the splitting $A=M-N$ with
\[
M=\frac{1}{\omega}(I-\gamma L)\quad \text{and} \quad   N=\frac{1}{\omega}[(1-\omega)I+(\omega-\gamma)L+\omega U),
\]
is an M-splitting of $A$. On the other hand,  we have $\rho(M^{-1}N)=\rho(\mathcal{L}_{\gamma,\omega})<1$. Therefore, by Lemma 1.1, we deduce that $A$ is an M-matrix. Now, the result follows immediately by Theorems 2.6 and 2.7 in \cite{Wang}. \quad $\Box$

In the sequel, we show that for improving  the convergence rate of the AOR iterative method, the preconditioner $\hat{P}=I+L+U$ is the best one between the preconditioners of the form $\tilde{P}=I+L(\alpha)+U(\alpha)$ with $\alpha_{ij}\in [0,1]$. We would like to point out here that if we set $\alpha_{ij}=1$  ($1\leq i \neq j \leq n$), then the preconditioner $\tilde{P}$ results in the preconditioner $\hat{P}$. Consider the AOR iteration matrix of the preconditioned system $\hat{P}Ax=\hat{P}b$ as follows:
\[
\hat{\mathcal{L}}_{\gamma,\omega} =(\hat{D}-\gamma\hat{L} )^{-1}[(1-\omega)\hat{D}+
(\omega-\gamma)\hat{L}+\omega\hat{U}],
\]
where $\hat{A}=\hat{P}A=\hat{D}-\hat{L}-\hat{U}$ in which  $\hat{D}$, $\hat{L}$ and $\hat{U}$ are the diagonal, strictly lower and strictly upper triangular matrices, respectively.

 {In the following, the set of indices $(i,j)$ associated with the nonzero off-diagonal entries of the matrix $A$ is represented by $\mathcal{N}_z(A)$, i.e.,
\[\mathcal{N}_z(A)=\{(i,j)| \quad i\neq j \quad \text{and} \quad a_{ij}\neq 0\}.\]}

\bigskip

\noindent\textbf{Theorem  2.3.} \textit{Suppose that $A=(a_{ij})\in \mathbb{R}^{n\times n}$ is a nonsingular Z-matrix.
Moreover, assume that $0\leq\gamma\leq\omega\leq 1$,  $\omega\neq 0$ and $\alpha_{ij} \in [0,1]$ for $i,j=1,2,\ldots,n$ ($i \neq j$). If
\begin{equation}\label{assump}
\sum_{k=1,k\neq i}^{n}(\alpha_{ik}-1) a_{ik}a_{kj}\leq 0, \quad  \text{for} \quad j < i \quad  {\text{and} \quad (i,j)\in \mathcal{N}_z(A)},
\end{equation}
and $\rho(\mathcal{L}_{\gamma,\omega})<1$, then
\begin{equation}\label{main}
\rho(\hat{\mathcal{L}}_{\gamma,\omega}) \leq \rho(\tilde{\mathcal{L}}_{\gamma,\omega}).
\end{equation}}
\noindent\textbf{Proof.} Let us consider the following splittings
 \[
 A=M-N, \qquad \tilde{A}=\tilde{M}-\tilde{N},  \qquad  \hat{A} = \hat{M} - \hat{N},
 \]
 where
\begin{eqnarray*}
  M         \hspace{-0.2cm}&=&\hspace{-0.2cm} \frac{1}{\omega}(I-\gamma L),\\
  N         \hspace{-0.2cm}&=&\hspace{-0.2cm} \frac{1}{\omega}[(1-\omega)I+(\omega-\gamma )L+\omega U],\\
  \tilde{M} \hspace{-0.2cm}&=&\hspace{-0.2cm} \frac{1}{\omega}(\tilde{D}-\gamma \tilde{L} ),\\
  \tilde{N} \hspace{-0.2cm}&=&\hspace{-0.2cm} \frac{1}{\omega} [(1-\omega)\tilde{D}+(\omega-\gamma )\tilde{L}+\omega\tilde{U}],\\
  \hat{M}   \hspace{-0.2cm}&=&\hspace{-0.2cm} \frac{1}{\omega} (\hat{D}-\gamma \hat{L} ),\\
  \hat{N}   \hspace{-0.2cm}&=&\hspace{-0.2cm} \frac{1}{\omega}[(1-\omega) \hat{D}+(\omega-\gamma )\hat{L}+\omega\hat{U}].
 \end{eqnarray*}
With an analogous strategy used in the proof of Theorem 2.2, we may show that the matrix $A$ is an M-matrix. Consequently, Theorem 2.1 shows that $\tilde{A}$ is also an M-matrix. Hence, the diagonal matrix $\tilde{D}$ has positive diagonal entries. Thus, the matrices $\tilde{M}$ and $\hat{M}$ are nonsingular and the splittings $\hat{A}=\hat{M}-\hat{N}$ and $\tilde{A}=\tilde{M}-\tilde{N}$ are M-splitting.

Considering the structures of $\tilde{A}$ and $\hat{A}$, we derive
\[
\tilde{D} - \hat{D}= E_1(1)-E_1(\alpha)  + E_2(1)-E_2(\alpha)  \geq 0,
\]
which is equivalent to say that $\tilde{D}\geq \hat{D}$.

For  {$(i,j)\notin \mathcal{N}_z(A)$}, we immediately conclude that
\[(\tilde{L} - \hat{L})_{ij} = \sum_{k=1,k\neq i,j}^{n}(\alpha_{ik}-1) a_{ik}a_{kj}\leq 0.\]
On the other hand, if  {$(i,j)\in \mathcal{N}_z(A)$},  Eq. \eqref{assump} implies that
\begin{eqnarray*}
(\tilde{L} - \hat{L})_{ij}=\sum_{k=1,k\neq i}^{n}(\alpha_{ik}-1) a_{ik}a_{kj}\leq 0.
\end{eqnarray*}
Thence, $(\tilde{L} - \hat{L})_{ij}\leq 0$ for $i=2,3,\ldots,n$ and $j=1,2,\ldots,i-1$ which means that $\tilde{L}\leq \hat{L}$.

As  $\rho(\gamma\tilde{D}^{-1}\tilde{L})< 1$, we may deduce that
\[
( \tilde{D}-\gamma\tilde{L})^{-1}=( I-\gamma\tilde{D}^{-1}\tilde{L})^{-1}\tilde{D}^{-1}
=I+\sum_{j=1}^{\infty} (\gamma\tilde{D}^{-1}\tilde{L})^j \tilde{D}^{-1}\geq 0.
\]
In a similar way, we can see that $(\hat{D}-\gamma\hat{L})^{-1}\geq 0$.

Straightforward computations show that $\hat{D}-\gamma\hat{L} \leq \tilde{D}-\gamma\tilde{L}$ which implies
\[
( \tilde{D}-\gamma\tilde{L})^{-1} \leq ( \hat{D}-\gamma\hat{L} )^{-1},
\]
or equivalently,
\begin{equation}\label{f1}
 0\leq \tilde{M}^{-1} \leq \hat{M}^{-1}.
\end{equation}
For the matrix $A$, we consider the following two splittings  $A=M_1-N_1=M_2-N_2$ where
\[
M_1 = \hat{P}^{ - 1} \hat{M}, \quad N_1 = \hat{P}^{ - 1} \hat{N}, \quad
 M_2= \tilde{P}^{ - 1} \tilde{M}  \quad \text{and} \quad N_2=\tilde{P}^{ - 1} \tilde{N}.
\]
By using Eq. \eqref{f1}, it can be verified that
\begin{equation}
M_1^{ - 1}=(\hat{P}^{ - 1} \hat{M})^{ - 1}  = \hat{M}^{ - 1} \hat{P} \ge \hat{M}^{ - 1} \tilde{P} \ge \tilde{M}^{ - 1} \tilde{P} = (\tilde{P}^{ - 1} \tilde{M})^{ - 1}=M_2^{ - 1}.
\end{equation}
It is not difficult to see that
\[\rho(M_1^{ - 1}N_1)=\rho(\hat{M}^{ - 1}\hat{N})<1 \quad  \text{and} \quad \rho(M_2^{ - 1}N_2)=\rho(\tilde{M}^{ - 1}\tilde{N})<1.\]
From Lemma 1.3, we deduce that
\[\rho(M_1^{ - 1}N_1)\leq \rho(M_2^{ - 1}N_2).\]
Or equivalently,
\[\rho(\hat{M}^{ - 1}\hat{N})\leq \rho(\tilde{M}^{ - 1}\tilde{N}),\]
which completes the proof. \quad $\Box$

We would like to comment here  that if $A$ is an irreducible matrix, then $A^{-1}>0$ \cite{Varga,Woznicki2} which implies strict inequality in \eqref{main}
when $ \tilde{M}^{-1} < \hat{M}^{-1}$, see Lemma 1.3.

\bigskip

\noindent\textbf{Theorem  2.4.} Let ${\mathcal{L}}_{\gamma,\omega}$ and $\tilde{\mathcal{L}}_{\gamma,\omega}$ denote the iteration matrices of the AOR and preconditioned AOR methods. Suppose that $A$ is an irreducible Z-matrix, $0\leq\gamma < 1$,  $\omega\neq 0$ and $\alpha_{ij} \in [0,1]$ for $i,j=1,2,\ldots,n$ ($i \neq j$). Moreover, assume that for each $(i,j)\in \mathcal{N}_z(A)$ there exists  $\tau\neq i,j$ such that
\[a_{ij}<\alpha_{ij}a_{ij}+\alpha_{i\tau}a_{i\tau}a_{\tau j} .\]
Then ${\mathcal{L}}_{\gamma,\omega}$  and $\tilde{\mathcal{L}}_{\gamma,\omega}$ are  nonnegative and irreducible matrices.

\bigskip

\noindent\textbf{Proof.} Let $A=I-L-U$ be an irreducible matrix, hence $\mathcal{G}(L+U)$ is strongly connected. We first show that $\tilde{A}=(\tilde{a}_{ij})$ is irreducible. To this end, using Lemma 1.4, we need to prove that $\mathcal{G}(\tilde{A})$ is strongly connected. Or equivalently, it is sufficient to demonstrate that  $\mathcal{N}_z(A)\subseteq \mathcal{N}_z(\tilde{A})$. For $(i,j)\in \mathcal{N}_z(A)$, by the assumption, we get
\[\tilde{a}_{ij}=a_{ij}-\displaystyle\sum_{k=1,k\neq i}^{n} \alpha_{ik}a_{ik}a_{kj}\leq a_{ij}-\alpha_{ij}a_{ij}- \alpha_{i\tau}a_{i\tau}a_{\tau j} < 0,\]
which is equivalent to say that $(i,j)\in \mathcal{N}_z(\tilde{A})$. Thus, the matrix $\tilde{A}$ is irreducible which implies that $\mathcal{G}(\tilde{D}^{-1}(\tilde{L}+\tilde{U}))$ is strongly connected where $\tilde{A}=\tilde{D}-\tilde{L}-\tilde{U}$.

It can be seen that if $A$ is a Z-matrix, then $\tilde{A}$ is a Z-matrix. Therefore, $\tilde{L},\tilde{U} \geq 0$ and we have:
\begin{eqnarray*}
\tilde{\mathcal{L}}_{\gamma,\omega}  \hspace{-0.2cm}&=&\hspace{-0.2cm}
(\tilde{D}-\gamma\tilde{L} )^{-1}[(1-\omega)\tilde{D}+(w-\gamma)\tilde{L}+\omega\tilde{U}] \\
                 \hspace{-0.2cm}&=&\hspace{-0.2cm}
(I-\gamma \tilde{D} ^{-1}\tilde{L} )^{-1}[(1-\omega)I+(\omega-\gamma) \tilde{D} ^{-1} \tilde{L}+\omega\tilde{D} ^{-1} \tilde{U}] \\
                 \hspace{-0.2cm}&=&\hspace{-0.2cm}  [\ I+(\gamma \tilde{D} ^{-1}\tilde{L} )+(\gamma \tilde{D} ^{-1}\tilde{L} )^{2} +\cdots \ ] [(1-\omega)I+(\omega-\gamma) \tilde{D} ^{-1} \tilde{L}+\omega\tilde{D} ^{-1} \tilde{U}]    \\
                    \hspace{-0.2cm}& \geq &\hspace{-0.2cm} [(1-\omega)I+(\omega-\gamma) \tilde{D} ^{-1} \tilde{L}+\omega\tilde{D} ^{-1} \tilde{U}]+(1-\omega)(\gamma \tilde{D} ^{-1}\tilde{L} ) \\
                    \hspace{-0.2cm}& = &\hspace{-0.2cm} [(1-\omega)I+ \omega(1-\gamma) \tilde{D} ^{-1} \tilde{L}+\omega\tilde{D} ^{-1} \tilde{U}].
\end{eqnarray*}
Hence, $\tilde{\mathcal{L}}_{\gamma,\omega}$ is a nonnegative and irreducible matrix.
In a similar manner, we may prove that
\[{\mathcal{L}}_{\gamma,\omega}\geq [(1-\omega)I+ \omega(1-\gamma) {L}+\omega {U}].\]
Using the facts that ${L},{U} \geq 0$ and $\mathcal{G}(L+U)$ is strongly connected, the above relation signifies that  ${\mathcal{L}}_{\gamma,\omega}$ is a nonnegative and irreducible matrix. \quad $\Box$

\section{Numerical experiments}
\label{SEC3}

All the numerical experiments presented in this section were
computed in double precision with some MATLAB codes on a Pentium 4 PC, with a 3.06 GHz CPU and
1.00GB of RAM.

\bigskip

\noindent \textbf{Example 1.} Consider the two dimensional convection-diffusion equation (see \cite{Wu})
\begin{equation}\label{eq005}
-(u_{xx}+u_{yy})+u_x+2u_y=f(x,y), \quad {\rm in} \quad \Omega = (0, 1) \times (0, 1),
\end{equation} with the homogeneous Dirichlet boundary conditions. Discretization of this equation on a uniform grid with $N \times N$ interior nodes ($n=N^2$), by using the second order centered differences for the
second and first order differentials, gives a linear system of equations of order $n$ with
$n$ unknowns. The coefficient matrix of the obtained linear system is of the form
\[
A=I \otimes P+Q \otimes I,
\]
where $\otimes$ denotes the Kronecker product,
\[
P=\textrm{tridiag}(-\frac{2+h}{8},1,-\frac{2-h}{8})\quad \textrm{and} \quad Q=\textrm{tridiag}(-\frac{1+h}{4},0,-\frac{1-h}{4}),
\]
are $N \times N$ tridiagonal matrices, in which the step size is $h=1/N$. We examine the following five preconditioners,
\begin{eqnarray*}
  P_0 \hspace{-0.2cm} &=& \hspace{-0.2cm} I,\\
  P_1 \hspace{-0.2cm} &=& \hspace{-0.2cm} I+0.5L, \\
  P_2 \hspace{-0.2cm} &=& \hspace{-0.2cm} I+0.5L+0.5U,\\
  P_3 \hspace{-0.2cm} &=& \hspace{-0.2cm} I+L(\alpha)+U(\alpha), \\
  {P_4} \hspace{-0.2cm} &{=}& \hspace{-0.2cm}
  {I+{\rm tridiag}(0,0,-a_{i,i+1}),\quad(\textrm{see~[4]})} \\
  {P_5} \hspace{-0.2cm} &{=}& \hspace{-0.2cm} {I+L+U},
\end{eqnarray*}
where for the preconditioner $P_3$, $\alpha_{ij}$'s are random numbers uniformly distributed in the interval $(0,1)$.
We would like to point out here that $P_0=I$ means that no preconditioner is used.

In \textsc{Figure} 1, we depict the eigenvalue distribution of $P_iA$  ($i=0,\ldots,5$) for $n=30^2=900$. This figure shows that the spectrum of the preconditioned matrix  $P_5A$ is more clustered than those of the matrices $P_iA$ for $i=0,1,2,3,4$.

In Table 1, the spectral radius of the AOR iterative method applied to the preconditioned systems $P_iAx=P_ib$ ($i=0,\ldots,5$) for different values of $\gamma$, $\omega$ and $n$ are given. As observed, the preconditioner  $P_5$ is the best one among the chosen  preconditioners.

For more investigation, we apply the GMRES($m$) method \cite{GMRES} with $m=10$ to solve $P_iAx=P_ib$ for $i=0,\ldots,5$.
In all of the experiments, the  vector $b=A(1,1,\ldots,1)^T$ was taken to be the right-hand side of the
linear system and a null vector as an initial guess. The stopping
criterion used was always
\[
\frac{\| b-Ax_k \|_2}{\| b \|_2} < 10^{-10}.
\]
In Table 2, we report the number of iterations and the CPU time (in parenthesis) for the convergence. As seen, the preconditioner $P_5$ is the best one.

\begin{figure}\label{Figure1}
\centering\includegraphics[height=5cm,width=7cm]{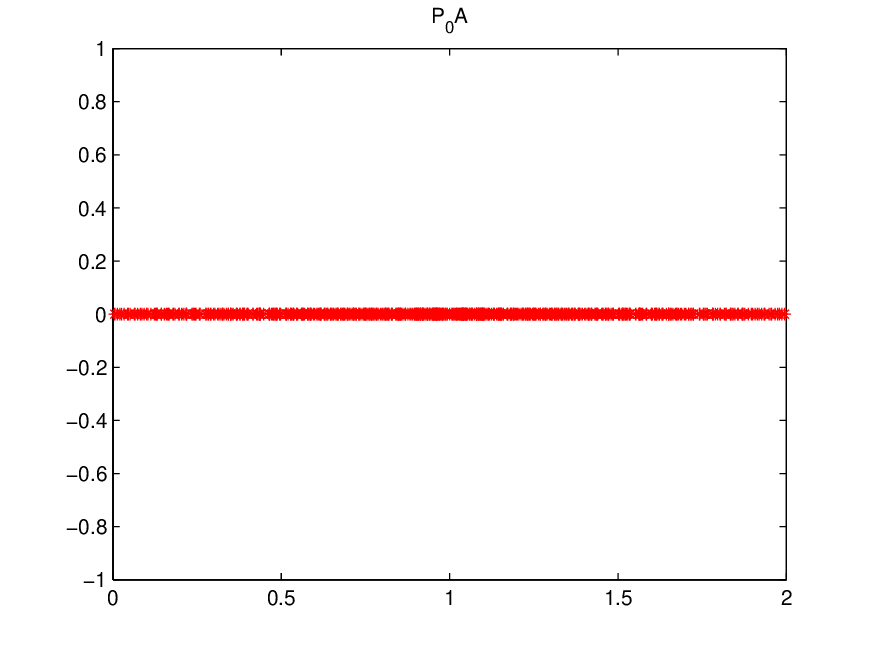}\includegraphics[height=5cm,width=7cm]{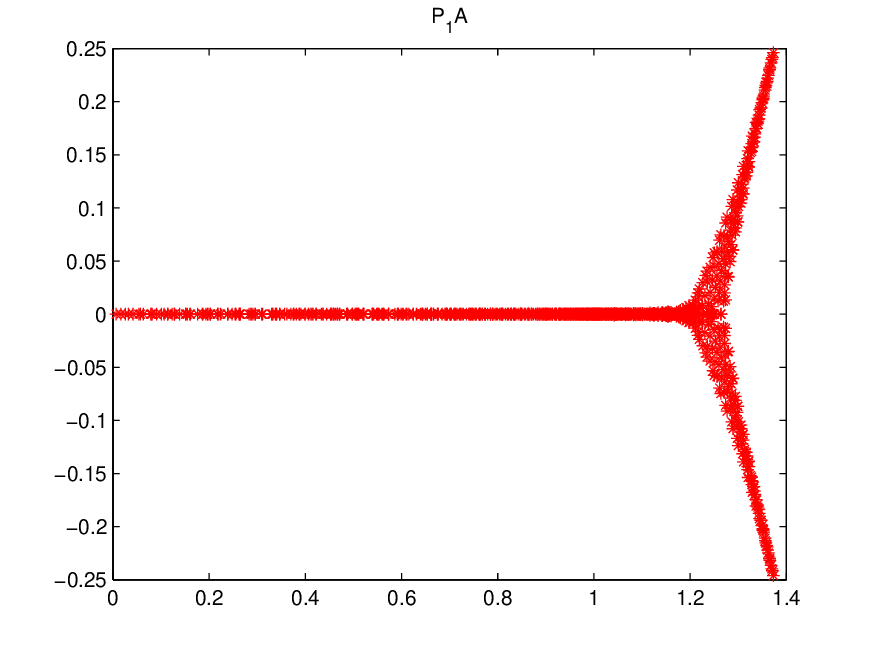}
\centering\includegraphics[height=5cm,width=7cm]{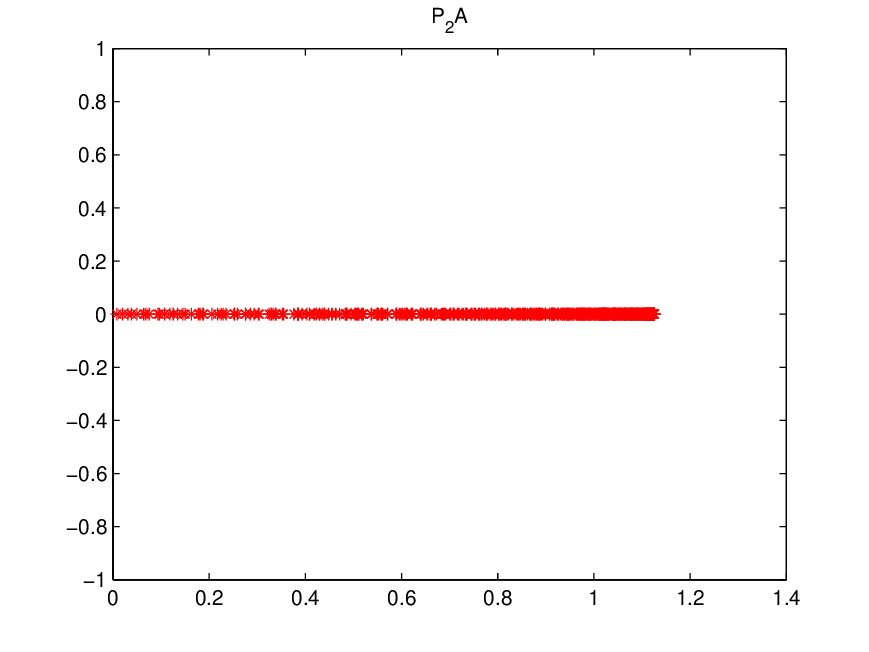}\includegraphics[height=5cm,width=7cm]{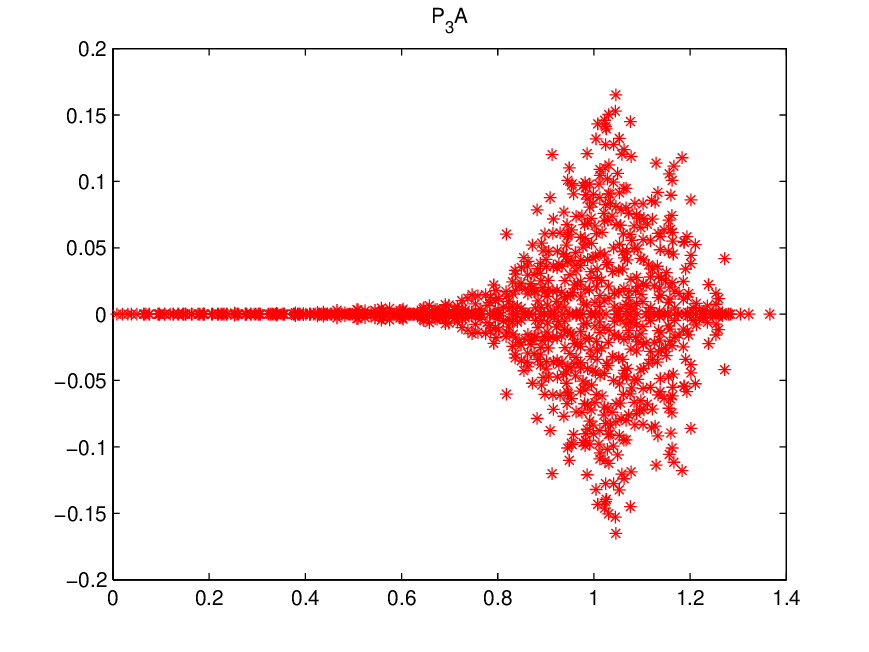}
\centering\includegraphics[height=5cm,width=7cm]{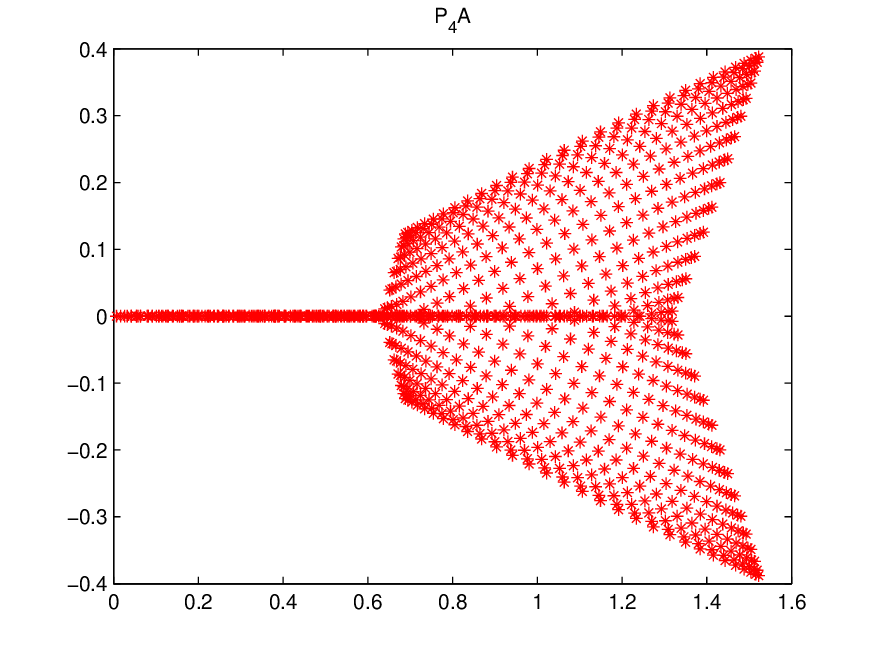}\includegraphics[height=5cm,width=7cm]{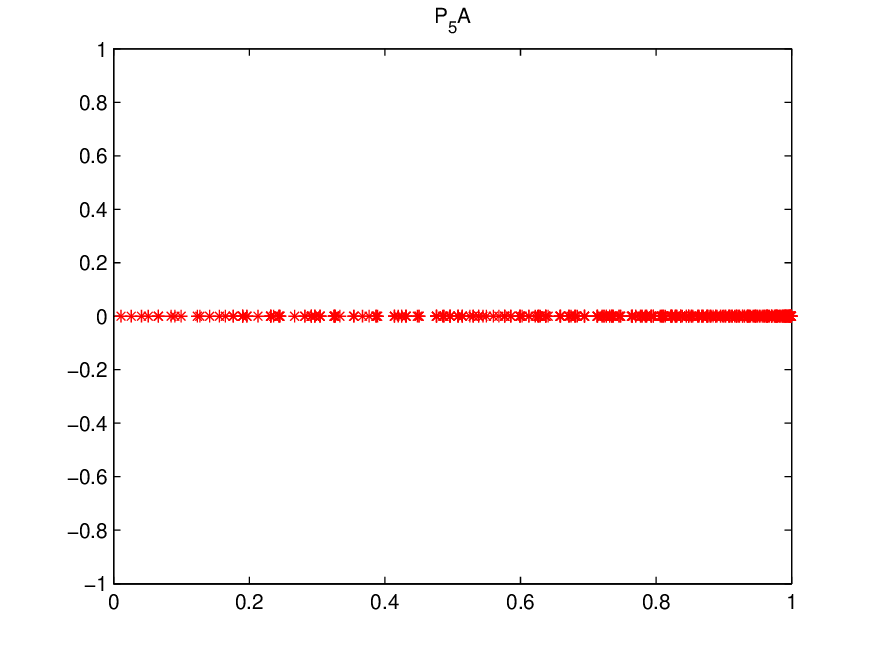}
\caption{Spectra of $P_iA$ for  Example 1 ($i=0,\ldots,5$). {Here, $P_0=I$ and} {$P_5=I+L+U$}. }
\bigskip
\bigskip
\end{figure}

\bigskip

\begin{table}\label{Table1}
\caption{Comparison of spectral radii for Example 1.   }
\vspace{-0.2cm}
\begin{center}
\begin{tabular}{lccccccc}\hline
$n~(=N^2)$  &  $(\gamma,\omega)$ &  $P_0$   & $P_1$   &  $P_2$  &   $P_3$  & $P_4$ & $P_5$  \\\hline
{25}     &  $(0.7,0.8)$       &  0.8317  & 0.7964  & 0.7404  &   0.7486 & {0.7657}  & 0.6323 \\
{25}    &  $(0.8,1)$         &  0.7739  & 0.7305  & 0.6540  &   0.6444 & {0.6798}  & 0.5138 \\
{100}   &  $(0.7,0.8)$       &  0.9474  & 0.9350  & 0.9125  &   0.9116 & {0.9230}  & 0.8677 \\
{100}   &  $(0.8,1)$         &  0.9289  & 0.9135  & 0.8821  &   0.8815 & {0.8933}  & 0.8221 \\\hline
\end{tabular}
\end{center}
\end{table}

\begin{table}\label{Tab2}
\caption{Number of iterations and the CPU time for the convergence of the GMRES(10) for Example 1.   }
\vspace{-0.2cm}
\begin{center}
\begin{tabular}{lcccccc}\hline
$n~(=N^2)$  &  $P_0$      & $P_1$      &   $P_2$     &    $P_3$   &    $P_4$   & $P_5$  \\\hline
{2500}   &  80(0.34)   & 57 (0.31)  &   33(0.22)  &  44(0.74)  & {79 (0.41)}  & 29(0.17) \\
{10000}  &  326(5.33)  & 130(2.83)  &  132(3.30)  &  110(2.75) & {191(3.88)}  & 78 (1.97)   \\
{22500}  &  702(29.08) & 365(19.73) &  244(15.56) &  350(23.02)& {534(31.36)} & 185(12.08) \\\hline
\end{tabular}
\end{center}
\end{table}

\noindent \textbf{Example 2.}  We consider the previous example with
\[
-(u_{xx}+u_{yy})+2e^{x+y}(xu_x+yu_y)=f(x,y), \quad {\rm in} \quad \Omega = (0, 1) \times (0, 1).
\]
All of the assumptions are the same as the previous example.

In \textsc{Figure} 2, we represent the eigenvalue distribution of $P_iA$ ($i=0,\ldots,5$) for $n=30^2=900$. This figure illustrates  that the spectrum of the preconditioned matrix  $P_5A$ is more clustered than those of the matrices $P_iA$ for $i=0,\ldots,4$.

In Table 3, the spectral radii of the AOR iterative method and in Table 4 numerical results of the GMRES(10) method  applied to the preconditioned systems  $P_iAx=P_ib$, $i=0,\ldots,5$ are given. As observed, the preconditioner  $P_5$ is the best one among the chosen  preconditioners.

\begin{figure}\label{Figure2}
\centering\includegraphics[height=5cm,width=7cm]{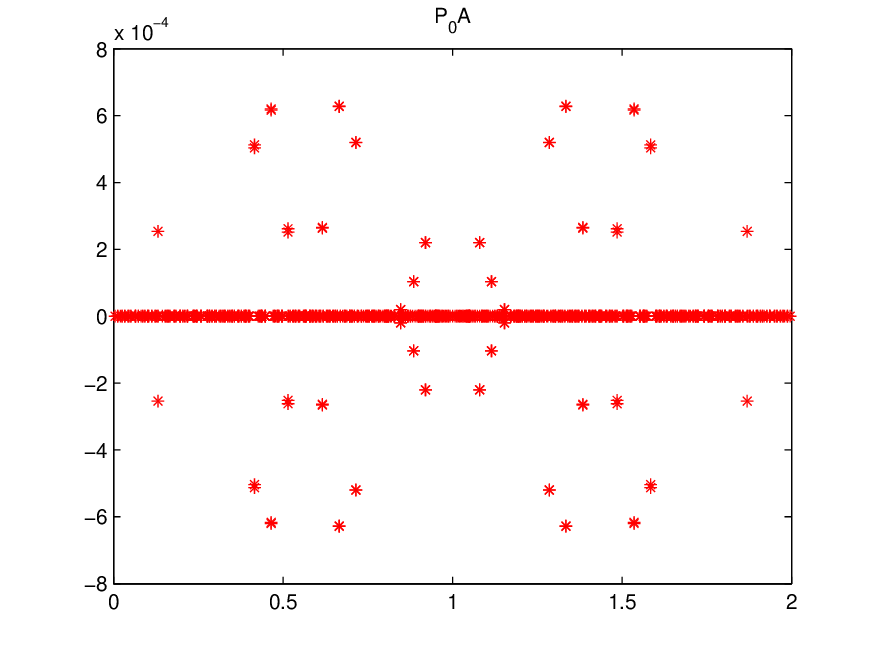}\includegraphics[height=5cm,width=7cm]{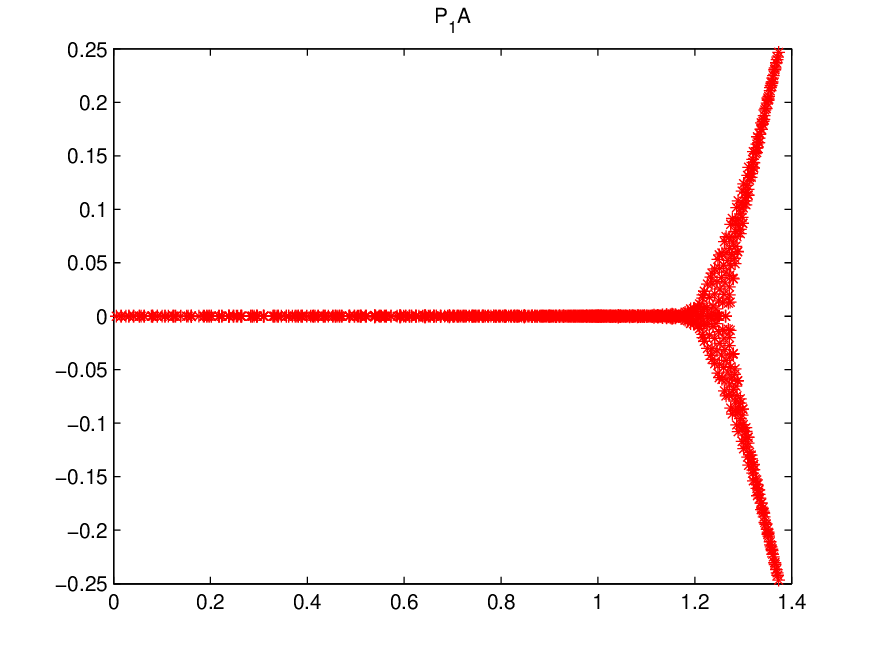}
\centering\includegraphics[height=5cm,width=7cm]{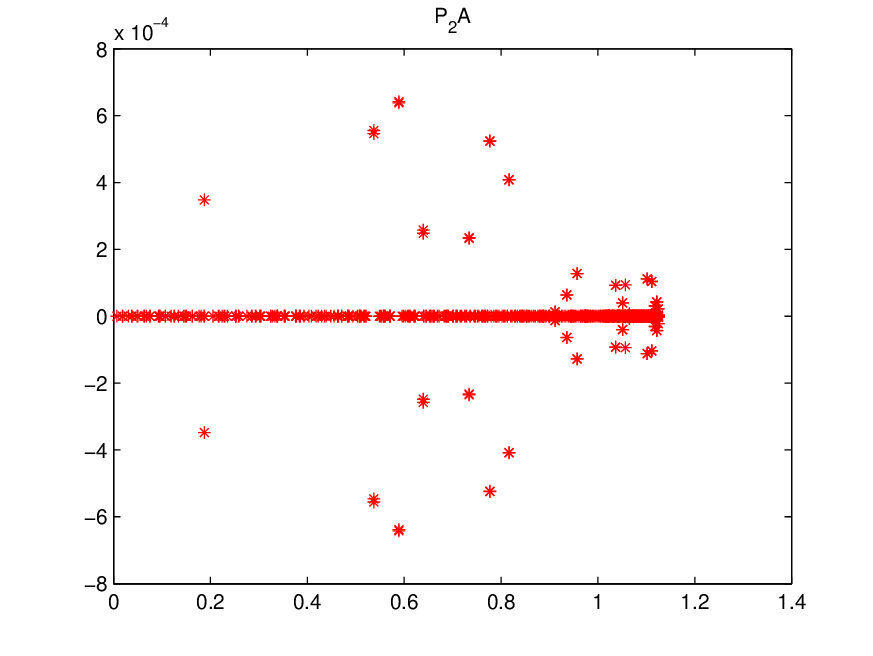}\includegraphics[height=5cm,width=7cm]{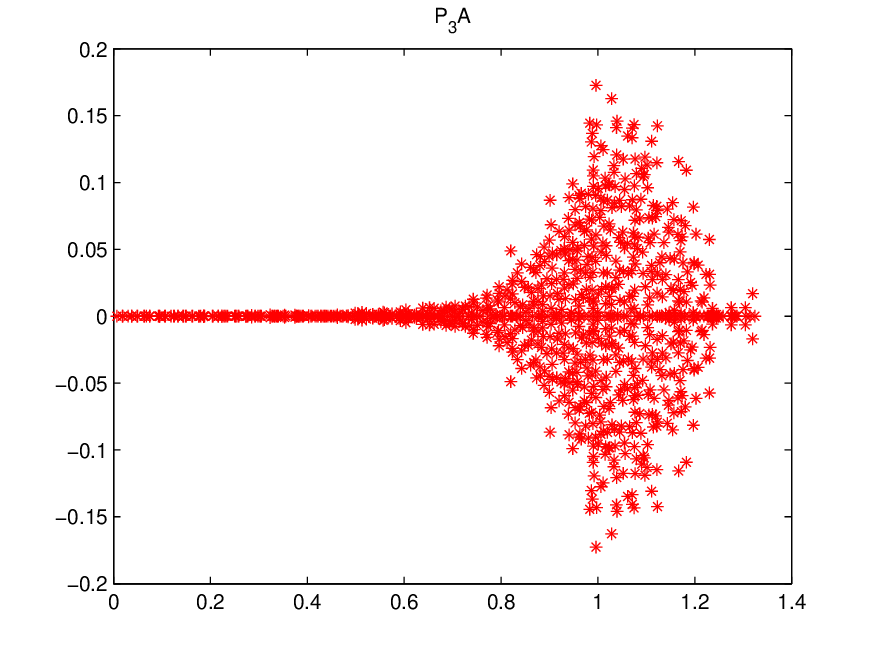}
\centering\includegraphics[height=5cm,width=7cm]{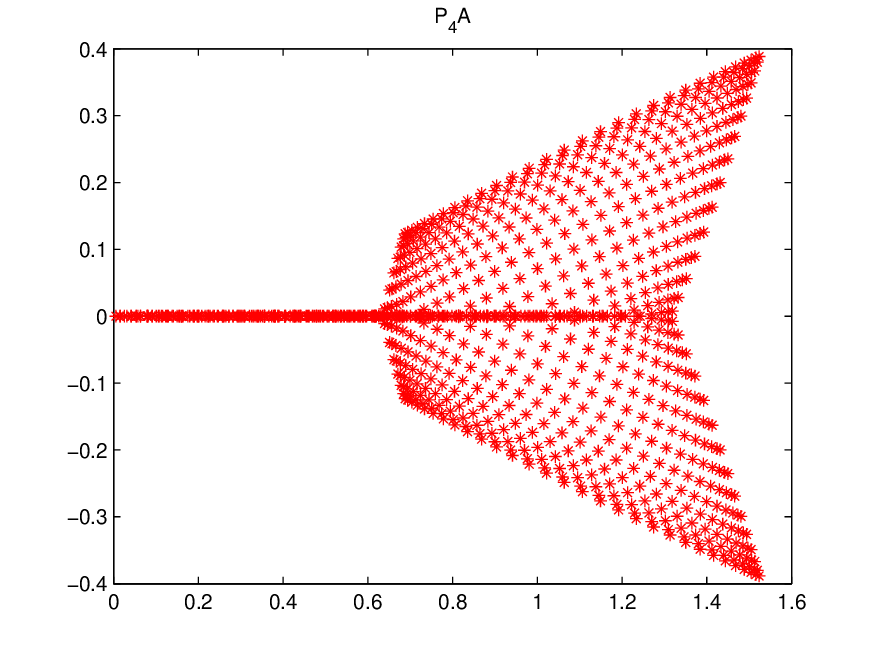}\includegraphics[height=5cm,width=7cm]{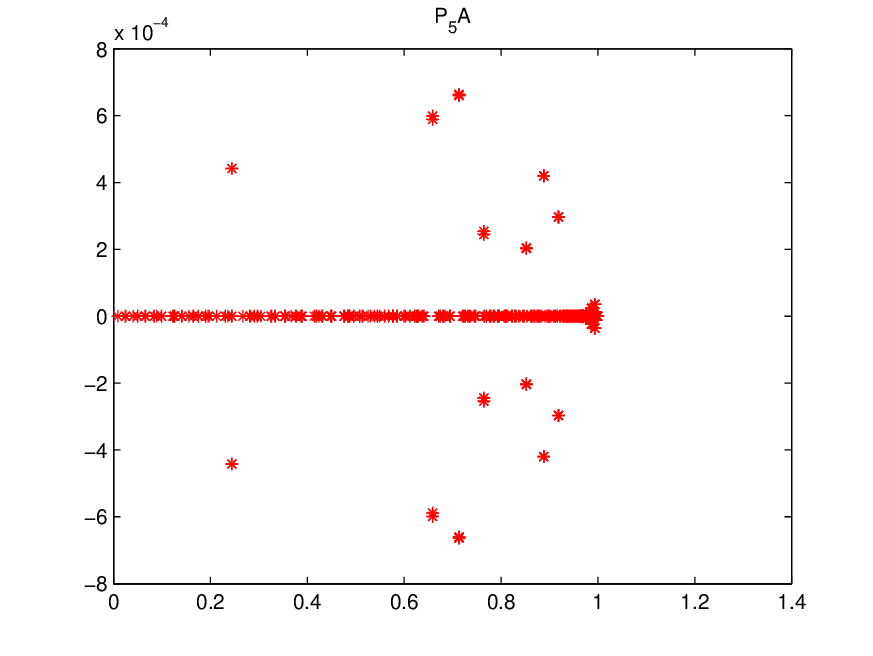}
\caption{Spectra of $P_iA$ for  Example 2 ($i=0,\ldots,5$). {Here, $P_0=I$ and} {$P_5=I+L+U$.} }
\bigskip
\bigskip
\end{figure}

\bigskip

\begin{table}\label{Tab3}
\caption{Comparison of spectral radii for Example 2.   }
\vspace{-0.2cm}
\begin{center}
\begin{tabular}{lccccccc}\hline
$n~(=N^2)$  &  $(\gamma,\omega)$ &  $P_0$  & $P_1$  &  $P_2$  &   $P_3$  & $P_4$& $P_5$  \\\hline
{25}    &  $(0.7,0.8)$       &  0.8657 & 0.8358 & 0.7871  &   0.8049 &  {0.8098}    & 0.6907 \\
{25}    &  $(0.8,1)$         &  0.8193 & 0.7823 & 0.7154  &   0.7033 &  {0.7390}    & 0.5891 \\
{100}   &  $(0.7,0.8)$       &  0.9581 & 0.9481 & 0.9298  &   0.9306 &  {0.9384}    & 0.8929 \\
{100}   &  $(0.8,1)$         &  0.9434 & 0.9309 & 0.9053  &   0.9045 &  {0.9145}    & 0.8558 \\\hline
\end{tabular}
\end{center}
\end{table}

\begin{table}\label{Tab4}
\caption{Number of iterations and the CPU time for the convergence of the
         GMRES(10) for Example 2.}
\vspace{-0.2cm}
\begin{center}
\begin{tabular}{lcccccc}\hline
$n~(=N^2)$  &  $P_0$      & $P_1$      &   $P_2$     &    $P_3$   & $P_4$     & $P_5$  \\\hline
{2500}          &  57(0.27)   & 46 (0.27)  &   28(0.19)  &  36 (0.22) & {54 (0.38)} & 23(0.13) \\
{3600}          &  85(0.52)   & 57 (0.45)  &   34(0.31)  &  49 (0.44) & {66 (0.65)} & 29(0.25)   \\
{4900}          &  92(0.72)   & 79(0.89)   &   49(0.61)  &  54 (0.72) & {101 (1.27)}& 37(0.47) \\
{6400}          &  111(1.23)  & 84(1.23)   &   52(0.84)  &  94 (1.55) & {103 (1.84)}& 45(0.73) \\\hline
\end{tabular}
\end{center}
\end{table}

\section{Conclusion and future work}

We have indicated that for improving the convergence rate of the AOR iterative method, the preconditioner $\hat{P}=I+L+U$ outperforms other preconditioners of the form $\tilde{P}=I+L(\alpha)+U(\alpha)$ with $\alpha_{ij}\in [0,1]$. Numerical experiments for the AOR and GMRES($m$) methods, to different preconditioned systems, have been reported to certify the established theoretical results.

A class of multi-level preconditioners and their associated preconditioned block AOR iterative method have been studied in \cite{Dou}. It can be easily verified that each of the preconditioners, exploited in a specific level of preconditioning,  pertains to a family of preconditioners which can be considered as a block form of our mentioned preconditioner. Future work may focus on comparing the multi-level preconditioned block AOR method with $\alpha_{ij}=1$ with the case that $\alpha_{ij}\neq 1$.

\end{document}